\documentclass[11pt,reqno,oneside]{amsart}
 \usepackage{latexsym}
 \usepackage{amssymb}
 \usepackage{amsfonts}
\usepackage{graphics}

\newtheorem{ansatz}{Ansatz}[section]
\newtheorem{conj}{Conjecture}[section]

\newtheorem{lem}{Lemma}[section]
\newtheorem{cor}{Corollary}[section]
\newtheorem{prop}{Proposition}[section]
\newtheorem{defi}{Definition}[section]
\newcommand{\R}{\mathbb{R}}
\newcommand{\N}{\mathbb{N}}

\renewcommand{\ker}{ker\,i(\alpha)}
\newcommand{\sd}{\mathcal{S}_{\delta}}
\newcommand{\s}{\mathcal{S}}

\newcommand{\cont}{{\mathcal{C}(M)}}
\newcommand{\fl}[2]{{\mathcal F}^{#1,#2}}
\newcommand{\gr}[2]{{\mathcal G}^{#1,#2}}

\newcommand{\dt}{\partial_{t}}

\title[Decomposition of symmetric tensor fields]{Decomposition of symmetric tensor fields in the presence of a flat contact projective structure}
 \author{Y. Fr\'egier, P. Mathonet, N. Poncin}
\date{\today}
\begin{document}
\begin{abstract} 
Let $M$ be an odd-dimensional Euclidean space endowed with a contact 1-form $\alpha$. We investigate the space of symmetric contravariant tensor fields over $M$ as a module over the Lie algebra of contact vector fields, i.e. over the Lie subalgebra made up of those vector fields that preserve the contact structure. If we consider symmetric tensor fields with coefficients in tensor densities (also called symbols), the vertical cotangent lift of the contact form $\alpha$ is a contact invariant operator. We also extend the classical contact Hamiltonian to the space of symmetric density valued tensor fields. This generalized Hamiltonian operator on the space of symbols is invariant with respect to the action of the projective contact algebra $sp(2n+2)$. These two operators lead to a decomposition of the space of symbols (except for some critical density weights), which generalizes a splitting proposed by V. Ovsienko in \cite{OC}.
\end{abstract}
\maketitle
\section{Introduction}
In a paper of 1997, C. Duval and V. Ovsienko \cite{DO2} considered the space ${\mathcal D}_{\lambda}(M)$ of differential operators acting on $\lambda$-densities ($\lambda\in\R$) on a manifold $M$ as modules over the Lie algebra of vector fields $\mathrm{Vect}(M)$. The density weight $\lambda$ allows to define a one parameter family of modules. The space of differential operators acting on half densities, which is very popular in mathematical physics in the context of geometric quantization, lies inside this family of representations. C. Duval and V. Ovsienko provided a first classification result for differential operators of order less or equal to 2. Differential operators that are allowed to modify the weight of their arguments also appear in the mathematical literature, for instance in projective differential geometry (see \cite{OvsBook0,Wil}). Hence the spaces ${\mathcal D}_{\lambda \mu}(M)$ (made of differential operators that map $\lambda$-densities into $\mu$-densities) also deserve interest. The first classification result of \cite{DO2} was then followed by a series of papers \cite{GO1,GO2,GA,LMT,MA} where the classification of the spaces ${\mathcal D}_{\lambda \mu}(M)$ was finally settled. More recently, in \cite{BHMP}, similar classification results were obtained for spaces of differential operators acting on differential forms.

The search for (local) $\mathrm{Vect}(M)$-isomorphisms from a space of differential operator to another one, or more generally the search of $\mathrm{Vect}(M)$-invariant maps between such spaces, can be made easier by using the so-called projectively equivariant symbol calculus introduced and developed by P. Lecomte and V. Ovsienko in \cite{LO} :
On the one hand the space ${\mathcal D}_{\lambda \mu}(\R^m)$ is naturally filtered by the order of differential operators. The action of $\mathrm{Vect}(\R^m)$ preserves the filtration and the associated graded space $\s_\delta(\R^m)$ (the space of symbols) is identified to contravariant symmetric tensor fields with coefficients in $\delta$- densities, with $\delta=\mu-\lambda$.

On the other hand the projective group $PGL(m+1,\R)$ acts locally on $\R^m$ by linear fractional transformations. The fundamental vector fields associated to this action generate a subalgebra of $\mathrm{Vect}(\R^m)$, the projective algebra $sl(m+1)$, which is obviously isomorphic to $sl(m+1,\R)$.

P. Lecomte and V. Ovsienko proved in \cite{LO} that the spaces ${\mathcal D}_\lambda(\R^m)$ and $\s_0(\R^m)$ are canonically isomorphic as $sl(m+1)$-modules. The canonical bijections are called respectively projectively equivariant quantization and symbol map.

This result was extended in \cite{Leclas} and \cite{DO1} : the spaces ${\mathcal D}_{\lambda \mu}(\R^m)$ and $\s_\delta(\R^m)$ are canonically isomorphic, provided $\delta$ does not belong to a set of critical values.

Once the projectively equivariant symbol map exists, the analysis of the filtered $sl(m+1)$-module ${\mathcal D}_{\lambda \mu}(\R^m)$ can be reduced to the study of the graded module $\s_\delta(\R^m)$.

Now, if we consider a contact manifold $M$ of dimension $m=2n+1$, it is sensible to consider the spaces ${\mathcal D}_{\lambda \mu}(M)$ as modules over the algebra of contact vector fields ${\mathcal C}(M)$, i.e. of those vector fields that preserve the contact structure. Using the projectively equivariant symbol map in order to perform local computations, we are lead to consider the space of symbols over $\R^{2n+1}$ as a representation of the ontact projective algebra ${\mathcal C}(\R^{2n+1})\cap sl(2n+2)= sp(2n+2)$ (this algebra is isomorphic to $sp(2n+2,\R)$).

A first question in the analysis of a representation is to know whether it decomposes as a direct sum of invariant subspaces or not.

The answer to this question is known for the space $\s^1_0$, that is, the space of vector fields : in  \cite{OC}, the author shows that the space of vector fields over a contact manifold, viewed as a module over the Lie algebra of contact vector fields, splits as the direct sum of contact vector fields and of vector fields which are tangent to the distribution :
\begin{equation}\label{splitval}\s^1_0(M)=\mathrm{Vect}(M)={\mathcal C}(M)\oplus TVect(M).\end{equation}
 The construction is based on the well-known Hamiltonian operator $X$, which associates a contact vector field to every $-\frac{1}{n+1}$-density, and is ${\mathcal C}(M)$-invariant.

In the present paper, we will define an extension of the operator $X$ to the whole space of symbols over the euclidean space $\R^{2n+1}$ endowed with its standard contact structure. This extended operator is not invariant with respect to the action of the algebra ${\mathcal C}(\R^{2n+1})$ but only with respect to the action of $sp(2n+2)$.

However, this operator, together with the contact form viewed as an operator of order zero acting on symbols, allows to define a decomposition of $\s_\delta(\R^{2n+1})$ as a sum of sp(2n+2)-submodules, unless $\delta$ belongs to a set of singular values.

The paper is organised as follows. In section \ref{sec2} we recall the definition of the basic material concerning densities, symbols and contact structures. We also recall the definition of the Lagrange bracket and of the operator $X$ acting on densities.

In section \ref{sec3} we show how to view the contact form as an invariant operator $i(\alpha)$. We give the definition of the extended operator $X$ on the space of symbols and show its invariance with respect to the action of the contact projective algebra $sp(2n+2)$.

 We show in the next section that these operators allow to define a representation of the Lie algebra $sl(2,\R)$ on the space of symbols. 

In section \ref{sec5}, we show that the contact hamiltonian operator $X$ allows to define an $sp(2n+2)$-invariant right inverse of the operator $i(\alpha)$, except for some singular values of the density weight. This right inverse allows to show the existence of a decomposition of the space of symbols into $sp(2n+2)$-invariant subspaces.

In the final section, we take another point of view and show how to obtain the decomposition by considering the natural filtration of the space of symbols of a given degree associated to the operator $i(\alpha)$.
\section{Basic objects}\label{sec2}
In this section, we recall the definitions of the basic objects that we will use throughout the paper, and we set our notations. As we continue, we denote by $M$ a smooth connected, Hausdorff and second countable manifold of dimension $m$. 
\subsection{Tensor densities and symbols}
Let us denote by $\Delta^{\lambda}(M)\to M$ the line bundle of tensor
densities of weight $\lambda$ over $M$ and by $\mathcal{F}_{\lambda}(M)$ of smooth sections of this bundle, i.e. the space $\Gamma(\Delta^{\lambda}(M))$. The Lie algebra of vector fields $\mathrm{Vect}(M)$ acts on $\mathcal{F}_{\lambda}(M)$ in a natural way. In local coordinates, any element $F$ of $\mathcal{F}_{\lambda}(M)$ as a local expression
\[F(x)=f(x)|dx^1\wedge\cdots\wedge dx^m|^\lambda\]
and the Lie derivative of $F$ in the direction of a vector field $X=\sum_iX^i\frac{\partial}{\partial x^i}$ is given by
\begin{equation}\label{ldens}(L_XF)(x)=(\sum_iX^i\frac{\partial}{\partial x^i}f+\lambda(\sum_i\frac{\partial}{\partial x^i}X^i)f)|dx^1\wedge\cdots\wedge dx^m|^\lambda.\end{equation}
Note that, as a vector space, $\mathcal{F}_{\lambda}(M)$ can be identified with the space of smooth functions on $M$, and thus formula (\ref{ldens}) defines a one parameter family of deformations of the natural representation of $\mathrm{Vect}(M)$ on $C^{\infty}(M).$
\subsection{Symbols}
We call the \emph{symbol space of degree $k$} and denote by $\s^k_\delta(M)$ (or simply $\s^k_\delta$)
the space of contravariant symmetric tensor fields of degree $k$, with
coefficients in $\delta$-densities :
\[\s^k_\delta(M) = \Gamma(S^kTM\otimes \Delta^{\delta}(M)).\]
We also consider the whole symbol space
\[
\s_\delta(M)=\bigoplus_{k\geq 0}\s^k_\delta(M).
\]
As we continue, we will freely identify symbols with
functions on $T^*M$ that are polynomial along the fibre and we will denote
by $\xi$ their generic argument in the fibre of $T^*M$.

The action of the algebra $\mathrm{Vect}(M)$ on symbols is the natural extension of its action on densities (\ref{ldens}) and on symmetric tensor fields. Let us write it down in order to illustrate the identification of symbols and functions on $T^*M$ :
\begin{equation}\label{lsymb}L_Xu(x,\xi)=\sum_i X^i\partial_{x^i}u(x,\xi)+\delta(\sum_i \partial_{x^i}X^i)u(x,\xi)-\sum_{i,k=1}^m(\partial_{x^k}X^i)\xi_i\partial_{\xi_k}u(x,\xi).\end{equation}
The spaces of symbols appear in a series of recent papers concerning equivariant quantizations \cite{LO,DLO,IFFT,MR}. Therefore we will not discuss them in full detail and refer the reader to these works for more information.
\subsection{Contact manifolds}
Here we will recall some basic facts about contact manifolds that we will use throughout the paper. Even though these facts are well known and exposed in various references (see for instance \cite{Ar1,Blair,McDuff}), we gather them for the paper to be self-contained and in order to fix notation. We will also focus our attention to the Euclidean space endowed with its standard contact structure since we are only concerned with local phenomena on contact manifolds.
\begin{defi}
A contact manifold is a manifold $M$ of odd dimension $m=2n+1$ together with a distribution of hyperplanes in the tangent space that is maximally non integrable (the contact distribution).
\end{defi}
Locally, the distribution can be defined as the kernel of a 1-form $\alpha$ and the (maximal) non-integrability condition means that $\alpha\wedge (d\alpha)^n\not=0.$ Moreover, all contact manifolds of dimension $m=2n+1$ are locally isomorphic to $\R^{2n+1}$ : there exist local coordinates (Darboux coordinates) such that the contact form writes
\begin{equation}\label{alpha}\alpha=\frac{1}{2}(\sum_{k=1}^n(p^kdq^k-q^kdp^k)-dt).\end{equation}
Unless otherwise stated, we will then only consider the Euclidean space $M=\R^{2n+1}$ with coordinates $(q^1,\ldots,q^n,p^1,\ldots,p^n,t)$ endowed with the contact form $\alpha$ defined by (\ref{alpha}).
\subsubsection{Contact vector fields}
A contact vector field over $M$ is a vector field which preserves the contact structure. The set of such vector fields forms a subalgebra of $\mathrm{Vect}(M)$, denoted by ${\mathcal C}(M)$. In other words, we have
\begin{equation}\label{contdef}{\mathcal C}(M)=\{Z\in\mathit{Vect}(M):\exists f_Z\in \Gamma(M\times \R^*) : L_Z\alpha=f_Z\alpha\}.\end{equation}
\subsubsection{The Lagrange bracket and the operator $X$ on densities}\label{Lag}
For every contact manifold $M$ there exists a ${\mathcal C}(M)$-invariant bidifferential operator acting on tensor densities. This is the so-called Lagrange bracket
\[\{,\}_{\mathcal L} :{\mathcal F}_{\lambda}(M)\times {\mathcal F}_{\mu}(M)\to{\mathcal F}_{\lambda+\mu+\frac{1}{n+1}}(M),\]
which is given in Darboux coordinates by the following expression :
\[\begin{array}{lll}\{f,g\}_{\mathcal L}&=&\sum_{k=1}^n(\partial_{p^k}f\partial_{q^k}g-\partial_{q^k}f\partial_{p^k}g)-\partial_tfE_s.g+\partial_tgE_s.f\\
&&+2(n+1)(\lambda f\partial_tg-\mu g\partial_tf),\end{array}\]
(where $E_s$ stands for the operator $\sum_{k=1}^n(p^k\partial_{p^k}+q^k\partial_{q^k}$)) for every $f\in {\mathcal F}_{\lambda}(M)$ and $g\in{\mathcal F}_{\mu}(M)$.
Now, the bilinear operator $\{,\}_{\mathcal L}$ can be viewed as a linear operator from ${\mathcal F}_{\lambda}(M)$ to the space ${\mathcal D}^1_{\mu, \lambda+\mu+\frac{1}{n+1}}(M)$ made of differential operators of order less or equal to one that map $\mu$-densities into $\lambda+\mu+\frac{1}{n+1}$-densities. Namely,
\[\{,\}_{\mathcal L} : f\mapsto \{f,\cdot\}_{\mathcal L} : g\mapsto \{f,g\}_{\mathcal L}.\]
Since the Lagrange bracket is ${\mathcal C}(M)$-invariant, this correspondence is also a ${\mathcal C}(M)$-invariant operator from ${\mathcal F}_{\lambda}(M)$ to ${\mathcal D}_{\mu, \lambda+\mu+\frac{1}{n+1}}(M)$ (the later space is endowed with the Lie derivative given by the commutator).

Finally, we consider the principal operator, which associates to every differential operator its term of highest order. It is well known that it is a $\mathrm{Vect}(M)$-invariant operator
\[\sigma : {\mathcal D}^1_{\mu, \lambda+\mu+\frac{1}{n+1}}(M)\to \s^1_{\lambda+\frac{1}{n+1}}(M).\]
If we compose these operators, we obtain
\[X:{\mathcal F}_{\lambda}(M)\to\s^1_{\lambda+\frac{1}{n+1}}(M): f\mapsto \sigma(\{f,\cdot\}_{\mathcal L}).\]
We then have the following immediate result.
\begin{prop}
The operator
\[X : {\mathcal F}_{\lambda}(M)\to\s^1_{\lambda+\frac{1}{n+1}}(M)\]
is $\cont$-invariant.
\end{prop}
Using the identification of symbols and polynomials, we can give the expression of $X$ :
\[X(f)(\xi)=\sum_{k=1}^n(\xi_{q^k}\partial_{p^k}f-\xi_{p^k}\partial_{q^k}f)-\partial_tf\langle E_s,\xi\rangle+\xi_t(E_s.f+2(n+1)\lambda f).\]
Let us close this section by the following result about contact vector fields.
\begin{prop}\label{contalg}
The algebra ${\mathcal C}(\R^{2n+1})$ is exactly $X({\mathcal F}_{\frac{-1}{n+1}}(\R^{2n+1})).$
\end{prop}
\subsubsection{The projective and symplectic algebras}
We also consider the projective Lie algebra $sl(2n+2)$. It is the
algebra of fundamental vector fields associated to the (local)
action of the projective group $PGL(2n+2,\R)$ on $\R^{2n+1}$. This
algebra is generated by constant and linear vector fields and
quadratic vector fields of the form $\eta{\mathcal E}$ for
$\eta\in\R^{2n+1^*}$, where ${\mathcal E}$ is the Euler vector field. Finally, the contact projective algebra
$sp(2n+2)$ is the intersection ${\mathcal C}(\R^{2n+1})\cap sl(2n+2)$. Note that the algebra $sp(2n+2)$ can be obtained using proposition \ref{contalg} by applying the operator $X$ to polynomial functions of degree less or equal to 2. For more details on the structure of this algebra, we refer the reader to \cite{Mat}.
\section{Invariant operators}\label{sec3}
Here we will define some operators related to the form $\alpha$. We
will then prove in the next sections that these operators commute
with the actions of ${\mathcal C}(M)$ or of $sp(2n+2)$.
\subsection{The contact form as an invariant operator}
We denote by $\Omega$ the volume form defined by $\alpha$, namely
\[\Omega= \alpha\wedge (d\alpha)^n\]
and by $div$ the divergence associated
$\Omega$ (actually the standard divergence over $\R^{2n+1}$). We then have for every $Z\in\cont$
\[L_Z\Omega=div(Z)\Omega=(n+1)f_Z\Omega,\]
(where $f_Z$ is defined in (\ref{contdef}))
and therefore
\[f_Z=\frac{1}{n+1}div (Z).\]
We then introduce a density weight in order to turn the form
$\alpha$ into a $\cont$- invariant tensor field : we consider
\[\alpha\otimes|\Omega|^{-\frac{1}{n+1}}\in\Gamma(T^*M\otimes \Delta^{-\frac{1}{n+1}}(M)).\]
and we can compute that the Lie derivative of this tensor field in
the direction of any field of $\cont$ is vanishing.

As we continue, we will omit the factor $|\Omega|^{-\frac{1}{n+1}}$
unless this leads to confusion, and consider $\alpha$ as an
invariant tensor field.

{\bf Remark :} This procedure is similar to the notion of conformal weight.

Now, since $\alpha$ is a 1-form, it defines a linear functional on vector
fields. This map has a natural extension to symmetric contravariant
tensor fields, which we denote by $i(\alpha)$. But since we want
this map to be $\cont$-invariant, we have to take the density weight of $\alpha$ into account and consider symmetric
tensor fields with coefficients in tensor densities.
We then have this first elementary result :
\begin{prop}\label{commutealpha}
The map
\[i(\alpha) : \sd\to\s^{k-1}_{\delta-\frac{1}{n+1}}\]
commutes with the action of $\cont$.
\end{prop}
We can give the expression of
the operator $i(\alpha)$ in terms of polynomials~: for every $S\in\sd^k$, there holds
\[i(\alpha)(S)(\xi)=\frac{1}{2}(\sum_i(p^i\partial_{\xi_{q^i}}-q^i\partial_{\xi_{p^i}})-\partial_{\xi_t})(i(\alpha)(S)).\]
{\bf Remarks :} The expression of the operator $i(\alpha)$ is
independent of the density weight. This weight appears only in order
to turn the map $i(\alpha)$ into an invariant map.
\subsection{The Hamiltonian operator $X$}
It turns out that the operator $X$ given in section \ref{Lag} extends to an operator on the space of symbols on $\R^{2n+1}$. We will prove that this operator is $sp(2n+2)$-invariant but not ${\mathcal C}(\R^{2n+1})$ equivariant.
\begin{defi}
We define the Hamiltonian operator $X$ in coordinates by
\[X:\sd^k\to\s^{k+1}_{\delta+\frac{1}{n+1}} : S\mapsto D(S)+a(k,\delta)\xi_t\,S,\]
where
\[D(S)(\xi)=\sum_i(\xi_{q^i}\partial_{p^i}- \xi_{p^i}\partial_{q^i})S(\xi)+\xi_tE_s(S)(\xi)-\langle E_s,\xi\rangle \partial_tS(\xi),\]
\[E_s=\sum_i(p^i\partial_{p^i}+q^i\partial_{q^i}),\]
and
\[a(k,\delta)=2(n+1)\delta-k.\]
\end{defi}
The main result is the following :
\begin{prop}
The operator $X : \sd^k\to\s^{k+1}_{\delta+\frac{1}{n+1}}$ commutes with the action of the algebra $sp(2n+2)$. It does not commute with the action of ${\mathcal C}(\R^{2n+1})$ unless $k=0$.
\end{prop}
\begin{proof}
The result could be checked by hands, and the computations are made easier by using the structure of the algebra $sp(2n+2)$, such as the grading of this algebra (see \cite{Mat}). However, we will present a few arguments based on the affine symbol map $\sigma_{\mathrm{Aff}}$ already used in \cite{LMT,IFFT,BHMP}. Recall that this map transforms a differential operator into a symmetric contravariant tensor field, just as the principal symbol map does. The map $\sigma_\mathrm{{Aff}}$ is equivariant with respect to the actions on differential operators and on symmetric tensor fields of the affine algebra $\mathrm{Aff}$, spanned by linear and constant vector fields. Hence the differential operator $X$ is $sp(2n+2)\cap\mathrm{Aff}$-invariant iff its affine symbol $\sigma_{\mathrm{Aff}}(X)$ is. 
We now compute 
\[(\sigma_{\mathrm{Aff}}(X)(\eta,S))(\xi)=(L_1(\xi,\eta) + a(k,\delta)\xi_t))S,\]
where
\[L_1(\xi,\eta)=\sum_i(\xi_{q^i}\eta_{p^i}- \xi_{p^i}\eta_{q^i})+\xi_t\langle E_s, \eta\rangle-\langle E_s,\xi\rangle \eta_t.\]
Moreover, the polynomial $L_1(\xi,\eta)$ is the principal symbol of the Lagrange bracket and is therefore an invariant polynomial.
It was already mentioned in \cite[section 6]{Mat} that the polynomial $P(\xi)=\xi_t$ is $sp(2n+2)\cap\mathrm{Aff}$-invariant. Therefore the operator $X$ commutes with the action of the subalgebra
$sp(2n+2)\cap \mathrm{Aff}$. 

Now, it is easy to verify that the algebra $sp(2n+2)$ is spanned by $sp(2n+2)\cap \mathrm{Aff}$ and the field
\[-\frac{1}{2}X_{t^2}=t E_s+t^2\partial_t=t E ,\]
where $E$ is the Euler field defined by $E=E_s+t\partial_t$.
The action of this field  on $\sd^{k}$ is given by the operator
\[t E - \langle E,\xi\rangle\dt + a(k,\delta) t ,\]
and it is easy to check that it commutes with $X$.

For the second part of the result, we check in the same way that the
commutator of $X$ with the Lie derivative in the direction of the
vector field $X_{q_1^3}$ does not vanish unless $k=0$.
\end{proof}
\section{Representation of the algebra $sl(2,\R)$}\label{sec4}
In general, the space of symbols is defined as the graded space
\[\s_{\delta}=\oplus_{k\in\N}\sd^k.\]
The grading is natural if we see the space of symbols as the graded
space associated to the filtered space of differential operators
mapping $\lambda$-densities into $\lambda+\delta$-densities.

However, this grading is not suitable to deal with the operators
$\alpha$ and $X$ since they modify the density weight of their
arguments. It is then natural to define the graded space
\[R_{\delta}=\oplus_{k\in\N}R^k_\delta\]
where $R^k_\delta=\s^k_{\delta+\frac{k}{n+1}}$ (and thus
$\sd^k=R^k_{\delta-\frac{k}{n+1}}).$

As we continue we will omit the reference to $\delta$ and denote
$R^k$ instead of $R^k_\delta$.

The operators $X$ and $i(\alpha)$ act on $R_\delta$. We can define a
degree on $R_\delta$ in such a way that $i(\alpha)$ has weight -1
and $X$ has weight +1.
\subsection{The operator $H$ and the representation of $sl(2,\R)$}
We now investigate the relationships between $X$ and $i(\alpha)$. We first compute the restriction of $X$ to $R^k$ and obtain the following elementary result
\begin{prop}\label{XRk}
The restriction of the operator $X$ to $R^k$ is given by
\[X(S)(\xi)=D(S)(\xi)+[2(n+1)\delta+k]\xi_t S(\xi)\]
for all $S\in R^k$.
\end{prop}
In order to make the results of this section easier to state, we introduce a new operator $H$ on $R_\delta$ :
\begin{defi}\label{defiH}
The operator $H$ is defined by its restrictions to $R^k$ given by
\[H\vert_{R^k}=h_k\mbox{Id}=-((n+1)\delta+ k)\mbox{Id}.\]
\end{defi}
The main result of this section deals with the commutators of the operators $X,i(\alpha)$ and $H$ on $R_{\delta}$.
\begin{prop}\label{sl2}
The operators $H$, $i(\alpha)$ and $X$ define a representation of the algebra $sl(2,\R)$ on the space $R_\delta$. Specifically, the relations
 \[\left\{\begin{array}{lcl}
   [i(\alpha),X] &=&  H\\
   {[H,i(\alpha)]}&=&i(\alpha)\\ 
  {[H,X]}&=&-X
 \end{array}\right.\]
hold on $R_\delta$.
\end{prop}
\begin{proof}
The first relation is a simple computation : on $R^k_\delta$, in view of Proposition \ref{XRk}, the commutator under consideration equals
\[\begin{array}{l}i(\alpha)\circ (D+(2(n+1)\delta+k)\xi_t)-(D+(2(n+1)\delta+k-1)\xi_t)\circ i(\alpha)\\
=[i(\alpha),D]+(2(n+1)\delta+k)[i(\alpha),\xi_t]+\xi_t\circ i(\alpha)\\
=-\frac{1}{2}E_\xi  -\frac{1}{2}(2(n+1)\delta+k)\mbox{Id} =h_k \mbox{Id},
\end{array}\]
(where $E_\xi$ stands for the operator $\sum_{i=1}^{2n+1}\xi^i\partial_{\xi^i}$) since the commutators in the second line are given by
\[ [i(\alpha),D]=-\frac{1}{2}E_{\xi}-\xi_t\circ \circ i(\alpha)\]
and
\[ [i(\alpha),\xi_t]=-\frac{1}{2}\mbox{Id}.\]
For the second relation we have on $R^k$ :
\[\begin{array}{lll}
[H,i(\alpha)]&=&H\circ i(\alpha)-i(\alpha) \circ H\\
&=&(h_{k-1}-h_k)i(\alpha)\\
&=&i(\alpha).\end{array}\]
The third one is proved in the same way.
\end{proof}
{\bf Remark :} In this representation of $sl(2,\R)$, $H$ corresponds to the action of an element in a Cartan subalgebra. The operator $i(\alpha)$ can be thought as the action of a generator of positive root space. In
this setting, the elements of $\ker$ could be seen as highest weight
vectors. 

We end this section with a technical result.
\begin{defi}
We define
\[r(l,k)=-\frac{l}{2}(2(n+1)\delta+2k+l-1).\]
\end{defi}
We can now state the result.
\begin{prop}\label{xl}
On $R^k$, we have 
\[\left\{
\begin{array}{lll}i(\alpha)\circ X^l-X^l\circ i(\alpha)&=&r(l,k)X^{l-1},\\
X\circ i(\alpha)^l-i(\alpha)^l\circ X&=&-r(l,k-l+1)i(\alpha)^{l-1}
\end{array}\right.\]
\end{prop}
\begin{proof}
There holds
\[i(\alpha)\circ X^l- X^l\circ i(\alpha)=\sum_{r=0}^{l-1}X^r[i(\alpha),X]X^{l-r-1}=\sum_{r=0}^{l-1}X^rHX^{l-r-1}.\]
On $R^k$, this operator is
\[\sum_{r=0}^{l-1}h_{k+l-r-1}X^{l-1}.\]
We then compute
\[\begin{array}{lll}\sum_{r=0}^{l-1}h_{k+l-r-1}&=&-\sum_{r=0}^{l-1}((n+1)\delta+k+l-r-1)\\
&=&-l((n+1)\delta+k+l-1)+\frac{l(l-1)}{2}\\&=&-\frac{l}{2}(2(n+1)\delta+2k+2l-2-l+1),
\end{array}\]
hence the first result.

We proceed in the same way for the second part :
We simply write
\[\begin{array}{lll}i(\alpha)^l\circ X-X\circ i(\alpha)^l&=&\sum_{r=0}^{l-1}i(\alpha)^r [i(\alpha),X] i(\alpha)^{l-r-1}\\
&=&(\sum_{r=0}^{l-1}h_{k-l+r+1}) i(\alpha)^{l-1}.
\end{array}\]
The result follows easily by the definition of $H$ :
\[\begin{array}{lll}
\sum_{r=0}^{l-1}h_{k-l+r+1}&=&-\sum_{r=0}^{l-1}((n+1)\delta+k-l+r+1)\\
&=&-\sum_{r=0}^{l-1}((n+1)\delta+k-r)\\
&=&-\frac{l}{2}(2(n+1)\delta+2k-l+1).
\end{array}\]
\end{proof}


\section{Decomposition of the space of symmetric tensors}\label{sec5}
In this section, we will obtain a decomposition of the space $R^k$ by induction on $k$. The idea is that $i(\alpha)$ maps $R^k$ to $R^{k-1}$. We will prove that this map is onto, and that there exist $sp(2n+2)$-invariant projectors $p_k : R^k\to R^k\cap \ker$, if $\delta$ does not belong to a set of singular values. These two fact will allow us to obtain the decomposition result.
\subsection{Expression of the projectors on $R^k\cap \ker$}\label{projsing}


The Ansatz for the shape of this projector is given by the
interpretation of $R_{\delta}$ as an $sl(2,R)$-module. More
precisely it comes from the following conjecture.

\begin{conj}
The algebra of all the invariant operators for the $sp(2n+2)$ action
on $R_\delta$ is generated by the operators $X$, $i(\alpha)$ and the projectors onto the spaces $R^k_\delta$.
\end{conj}

In particular if we want to build a projector $\Pi_k$ from $R^k$
onto a submodule of itself, in the category of modules it should be a linear combination of
compositions of $X$ and $i(\alpha)$. Moreover all the monomials in this expressions should have the same degree in $X$ and $i(\alpha)$. Finally, we can order these compositions, using Proposition \ref{xl} and the projector should be of the form $\Pi_k= \sum_{l=0}^{\infty}
b_{k,l}X^l\circ i(\alpha)^l $. But since $i(\alpha)^l_{\vert R^k}=0$ for $l>k$. Hence one gets the
following Ansatz :

\begin{ansatz}
The projector $p_k : R^k\to R^k\cap \ker $ should be of the form
\begin{equation}\label{ans}p_k=\mbox{Id}+\sum_{l=1}^kb_{k,l}X^l\circ i(\alpha)^l.\end{equation}
\end{ansatz}

It is possible to determine 
the expression of the constants $b_{k,l}$ from the constraint
 $p_k^2=p_k$. In the following we give directly their expression and check that
 we effectively obtain a projector. It turns out that there exists a set $I_k$ of values of $\delta$ such that
 no operator of the shape (\ref{ans}) can be a projector on $R^k\cap \ker$. In this situation, if the conjecture is true, the space $R^k$ is not the direct sum of $R^k\cap \ker$ and of an $sp(2n+2)$ invariant subspace.
\begin{defi}
For every $k\geq 1$, we set
 \[I_k=\{-\frac{p}{2(n+1)}:p\in\{k-1,\cdots,2k-2\}\}.\]
\end{defi}
We then have the following result.
\begin{prop}
If $\delta\not\in I_k$, the operator $p_k: R^k\to R^k$ defined by
\begin{equation}\label{proj}p_k :R^k\to R^k : p_k=\mbox{Id}+\sum_{l=1}^kb_{k,l}X^l\circ i(\alpha)^l\end{equation}
where
\[b_{k,l}=(\Pi_{j=1}^l-r(j,k-j))^{-1},\]
is a projector onto $R^k\cap \ker$.
\end{prop}
\begin{proof}
The restriction of $p_k$ to $\ker$ is the identity mapping, in view
of (\ref{proj}). It is then sufficient to prove that
$Im\,(p_k)\subset \ker$, and we then deduce $Im\,(p_k)=\ker$ and
then $p_k^2=p_k$.

We actually have the relation
\[i(\alpha)\circ p_k=0.\]
Indeed, we have
\[\begin{array}{lll}i(\alpha)\circ p_k&=&i(\alpha)+\sum_{l=1}^k b_{k,l}i(\alpha) X^l i(\alpha)^l\\
&=&i(\alpha)+\sum_{l=1}^k b_{k,l}(X^l i(\alpha)^{l+1}+r(l,k-l)X^{l-1} i(\alpha)^l).
\end{array}\]
The result follows since the constants $b_{k,l}$ fulfill the
equations
\[\left\{\begin{array}{lll}
1+b_{k,1}r(1,k-1)&=&0\\
b_{k,l}+b_{k,l+1}r(l+1,k-l-1)&=&0\quad\forall l=1,\cdots k-1,
\end{array}\right.\]
and since $i(\alpha)^{k+1}\equiv 0$ on $R^k$.
\end{proof}

This proposition suggests to define a new operator from $R^{k-1}$ to
$R^k$.
\begin{defi}
If $\delta\not\in I_k$, we define
\[s_{k-1}:R^{k-1}\to R^k : S\mapsto (-\sum_{l=1}^{k}b_{k,l}X^l i(\alpha)^{l-1})(S).\]
\end{defi}
We then have the following result that links $s_{k-1}$ and
$i(\alpha)$ :
\begin{lem}\label{inv}
If $\delta\not\in I_k$, then one has $i(\alpha)\circ s_{k-1}=\mbox{Id}$ on $R^{k-1}$.
\end{lem}
\begin{proof}
We just compute
\[\begin{array}{lll}
i(\alpha)\circ s_{k-1}&=&-\sum_{l=1}^{k}b_{k,l}i(\alpha) X^li(\alpha)^{l-1}\\
&=&-\sum_{l=1}^{k}b_{k,l}(X^li(\alpha)^l+r(l,k-l)X^{l-1}i(\alpha)^{l-1})\\&=&\mbox{Id},
\end{array}\]
by using Proposition \ref{xl}, and the definition of $b_{k,l}$.
\end{proof}
From this Lemma, we obtain an important information about $i(\alpha)$.
\begin{cor}\label{surj}
For every $k$ and $\delta$, the map $i(\alpha) : R^k_\delta\to R^{k-1}_\delta$ is
surjective.
\end{cor}
\begin{proof}
If $\delta\not\in I_k$, the result follows from the existence of a
right inverse. Since the expression of $i(\alpha)$ is independent of
$\delta$, the result holds true for every $\delta$.
\end{proof}
We can now state a first decomposition result.
\begin{prop}
If $\delta\not\in I_k$, one has $R^k=\ker\oplus s_{k-1}(R^{k-1})$.
\end{prop}
\begin{proof}
Since the projector $p_k$ is defined, there is a decomposition
\[R^k=(\ker\cap R^k)\oplus V^k,\] where $V^k=Im(\mbox{Id}-p_k)$. It
follows that the restriction of $i(\alpha)$ to $V^k$ is injective. It
is also surjective by corollary \ref{surj}. By Lemma \ref{inv},
$s_{k-1}$ is the inverse of $i(\alpha)$ and thus
$V^k=s_{k-1}(R^{k-1})$.
\end{proof}
By applying successively the previous proposition we obtain a second
result. We denote by $s$ the operator on $R_\delta$ whose restriction to $R^{k-1}$ is $s_{k-1}$.
\begin{cor}
If $\delta\not\in \cup_{j=1}^kI_j,$ then there holds
\[R^k=\oplus_{l=0}^k s^l(R^{k-l}\cap \ker).\]
\end{cor}
We will now compute the restriction of $s^l$ to $(R^{k-l}\cap
\ker)$ and show that it is a scalar multiple of  $X^l$.
\begin{prop}
Suppose that $\delta\not\in \cup_{j=1}^kI_j$. Then the
restriction of $s^l$ to $R^{k-l}\cap \ker$ equals $c(l,k-l)
X^l$, where
\[c(l,k-l)=(\Pi_{i=1}^lr(i,k-l))^{-1}.\]
\end{prop}
\begin{proof}
We first prove the existence of the constant $c(l,k-l)$ by showing
that the restriction of $s^j$ to $R^{k-l}\cap \ker$ equals
$c(j)X^j$ $(c(j)\in\R)$, for all $j=1,\cdots, l$. For $j=1$, the
result follows from the very definition of $s_{k-l}$, and we
obtain $c(1)=-b_{k-l+1,1}$. Suppose that it holds true for $s^j$.
We then have on $R^{k-l}\cap \ker$
\[s^{j+1}=s\circ s^j=-\sum_{a=1}^{k-l+j+1}b_{k-l+j+1,a}X^ai(\alpha)^{a-1} s^j.\]
The last term is a multiple of $X^{j+1}$, by induction, since the
composition $\alpha^{a-1}s^j$ vanishes if $a-1>j$ and is equal
to $s^{j-a+1}$ if $a-1\leq j$.

 Now, we compute the value of $c(j)$ by analysing the restriction of the operator $i(\alpha)\circ s^{j+1}$ to $R^{k-l}\cap \ker$. On the one hand, it is equal to $s^j$, that is to $c(j)X^j$. On the other hand, it writes
\[i(\alpha)\circ s^{j+1}=c(j+1)i(\alpha) X^{j+1}=c(j+1)(X^{j+1}i(\alpha) +r(j+1,k-l)X^j),\]
by Proposition \ref{xl}.

Finally, we obtain the relation
\[(c(j)-r(j+1,k-l)c(j+1))X^j=0,\]
on $R^{k-l}\cap \ker$, and therefore, since $X^j$ is injective on
this space
\[c(j)=r(j+1,k-l)c(j+1),\]
and the result follows.
\end{proof}
\section{Another point of view}
In this section we investigate the relations of the maps $i(\alpha)$ and $X$ with the filtration induced by $i(\alpha)$. We denote by $\fl{k}{l}$ the space $R^k\cap ker\,i(\alpha)^l$. Since $i(\alpha)$ is a $\cont$-invariant operator, the spaces $\fl{k}{l}$ are stable under the action of $\cont$. 
Moreover there is an obvious filtration of $R^k$ defined by
\[0=\fl{k}{0}\subset \fl{k}{1}\subset\cdots\subset \fl{k}{k+1}=R^k.\]
We first prove that the maps $i(\alpha)$ and $X$ respect this filtration and therefore induce mappings that deserve special interest on the associated graded spaces.
\begin{prop}\label{ialphax}
We have 
\[\left\{
\begin{array}{lll}
i(\alpha)(\fl{k}{l})&\subset& \fl{k-1}{l-1}\\
X(\fl{k-1}{l-1})&\subset& \fl{k}{l}
\end{array}\right.
\]
for every $k\geq 1$ and every $l\in\{1,\cdots,k+1\}$.
\end{prop}
\begin{proof}
The first result is a direct consequence of the definition of the filtration. For the second one, suppose that $v$ is an element of $\fl{k-1}{l-1}$ and compute, using Proposition \ref{xl},
\[i(\alpha)^l\circ X(v)=X\circ i(\alpha)^l (v)+ r(l,k-l)i(\alpha)^{l-1}(v).\]
It follows that $X(v)$ is in $\fl{k}{l}$.
\end{proof}
Let us introduce some more notation.
\begin{defi}
For every $l\in\{1,\cdots,k+1\}$ we denote by $\gr{k}{l}$ the quotient space $\fl{k}{l}/ \fl{k}{l-1}$. This space is naturally endowed with a representation of $\cont$. In particular we have $\gr{k}{1}\cong\fl{k}{1}$. We also set $\gr{k}{0}=\{0\}$.
\end{defi}
By Proposition \ref{ialphax}, the maps $i(\alpha)$ and $X$ induce maps on the graded space $\oplus_{k,l}\gr{k}{l}$. Namely, for every $l\in\{1,\cdots,k+1\}$, we set :
\[\left\{\begin{array}{l}
\widetilde{i(\alpha)} : \gr{k}{l} \to \gr{k-1}{l-1}:[u]\mapsto [i(\alpha)(u)]\\
\widetilde{X}:\gr{k-1}{l-1}\to \gr{k}{l}:[u]\mapsto [X(u)]
\end{array}\right.\] 
The main property of these maps is the following.
\begin{prop}\label{inverse}
For every $k\geq 1$ and $l\in\{1,\cdots,k+1\}$, we have 
\[\widetilde{X}\circ\widetilde{i(\alpha)}\vert_{\gr{k}{l}}=r(l-1,k-l+1)Id,\]
and
\[\widetilde{i(\alpha)}\circ\widetilde{X}\vert_{\gr{k-1}{l-1}}=r(l-1,k-l+1)Id.\]
In particular, if $r(l-1,k-l+1)$ does not vanish, the restricted map
$\widetilde{i(\alpha)} : \gr{k}{l}\to \gr{k-1}{l-1}$ is invertible and the inverse map is given by
$\frac{1}{r(l-1,k-l+1)}\widetilde{X}.$
\end{prop}
\begin{proof}We only prove the first identity. The second one can be proved using the same arguments.

Let $u$ be in $\fl{k}{l}$. By definition, we have
\[\widetilde{X}\circ\widetilde{i(\alpha)}([u])=[X\circ i(\alpha)(u)].\]
By Proposition \ref{xl}, we have 
\[\begin{array}{lll}i(\alpha)^{l-1}\circ X\circ i(\alpha)(u)&=&X\circ i(\alpha)^l(u)+r(l-1,k-l+1)i(\alpha)^{l-1}(u)\\
&=&r(l-1,k-l+1)i(\alpha)^{l-1}(u),
\end{array}\]
so that $[X\circ i(\alpha)(u)]=r(l-1,k-l+1)[u]$, and the result follows.
\end{proof}
We get an immediate corollary about the decomposition of the filter $\fl{k}{l}$ into stable subspaces.
\begin{cor}
Consider $k\geq 1$ and $l\in\{2,\cdots,k+1\}$ and suppose that $\prod_{j=1}^{l-1}r(j,k-l+1)\not=~0$. There holds
\[\fl{k}{l}=\fl{k}{l-1}\oplus X^{l-1}(\fl{k-l+1}{1}).\]
\end{cor}
\begin{proof}
Consider the following short exact sequence of $sp(2n+2)$-modules
\[0\longrightarrow\fl{k}{l-1}\longrightarrow\fl{k}{l}\longrightarrow\gr{k}{l}\longrightarrow 0.\]
By the previous proposition the $sp(2n+2)$-modules $\gr{k}{l}$ and $\gr{k-l+1}{1}$ are isomorphic through $\widetilde{i(\alpha)}^{l-1}$ and $\widetilde{\phi}=\frac{1}{\prod_{j=1}^{l-1}r(j,k-l+1)}\widetilde{X}^{l-1}$. Denote by $t$ the trivial isomorphism $\gr{k-l+1}{1}\to \fl{k-l+1}{1}$. It is then easy to check that the map
\[\phi=\frac{1}{\prod_{j=1}^{l-1}r(j,k-l+1)}X^{l-1}\circ t\circ \widetilde{i(\alpha)}^{l-1}\]
provides a section of the exact sequence above.
\end{proof}

\def\cprime{$'$}


\begin{thebibliography}{10}

\bibitem{Ar1}
V.~I. Arnold.
\newblock {\em Mathematical methods of classical mechanics}.
\newblock Springer-Verlag, New York, 1978.
\newblock Translated from the Russian by K. Vogtmann and A. Weinstein, Graduate
  Texts in Mathematics, 60.

\bibitem{Blair}
David~E. Blair.
\newblock {\em Contact manifolds in {R}iemannian geometry}.
\newblock Springer-Verlag, Berlin, 1976.
\newblock Lecture Notes in Mathematics, Vol. 509.

\bibitem{BHMP}
F.~Boniver, S.~Hansoul, P.~Mathonet, and N.~Poncin.
\newblock Equivariant symbol calculus for differential operators acting on
  forms.
\newblock {\em Lett. Math. Phys.}, 62(3):219--232, 2002.

\bibitem{IFFT}
F.~Boniver and P.~Mathonet.
\newblock I{FFT}-equivariant quantizations.
\newblock {\em J. Geom. Phys.}, 56(4):712--730, 2006.

\bibitem{DLO}
C.~Duval, P.~Lecomte, and V.~Ovsienko.
\newblock Conformally equivariant quantization: existence and uniqueness.
\newblock {\em Ann. Inst. Fourier (Grenoble)}, 49(6):1999--2029, 1999.

\bibitem{DO1}
C.~Duval and V.~Ovsienko.
\newblock Projectively equivariant quantization and symbol calculus:
  noncommutative hypergeometric functions.
\newblock {\em Lett. Math. Phys.}, 57(1):61--67, 2001.

\bibitem{DO2}
C.~Duval and V.~Yu. Ovsienko.
\newblock Space of second-order linear differential operators as a module over
  the {L}ie algebra of vector fields.
\newblock {\em Adv. Math.}, 132(2):316--333, 1997.

\bibitem{GA}
H.~Gargoubi.
\newblock Sur la g\'eom\'etrie de l'espace des op\'erateurs diff\'erentiels
  lin\'eaires sur {$\bold R$}.
\newblock {\em Bull. Soc. Roy. Sci. Li\`ege}, 69(1):21--47, 2000.

\bibitem{GO2}
H.~Gargoubi and V.~Yu. Ovsienko.
\newblock Space of linear differential operators on the real line as a module
  over the {L}ie algebra of vector fields.
\newblock {\em Internat. Math. Res. Notices}, (5):235--251, 1996.

\bibitem{GO1}
Kh. Gargubi and V.~Ovsienko.
\newblock Modules of differential operators on the line.
\newblock {\em Funktsional. Anal. i Prilozhen.}, 35(1):16--22, 96, 2001.

\bibitem{LMT}
P.~B.~A. Lecomte, P.~Mathonet, and E.~Tousset.
\newblock Comparison of some modules of the {L}ie algebra of vector fields.
\newblock {\em Indag. Math. (N.S.)}, 7(4):461--471, 1996.

\bibitem{LO}
P.~B.~A. Lecomte and V.~Yu. Ovsienko.
\newblock Projectively equivariant symbol calculus.
\newblock {\em Lett. Math. Phys.}, 49(3):173--196, 1999.

\bibitem{Leclas}
Pierre B.~A. Lecomte.
\newblock On the projective classification of the modules of differential
  operators on {$\Bbb R\sp m$}.
\newblock In {\em Noncommutative differential geometry and its applications to
  physics (Shonan, 1999)}, volume~23 of {\em Math. Phys. Stud.}, pages
  123--129. Kluwer Acad. Publ., Dordrecht, 2001.

\bibitem{MA}
P.~Mathonet.
\newblock Intertwining operators between some spaces of differential operators
  on a manifold.
\newblock {\em Comm. Algebra}, 27(2):755--776, 1999.

\bibitem{Mat}
P.~Mathonet.
\newblock Invariant bidifferential operators on tensor densities over a contact
  manifold.
\newblock {\em Lett. Math. Phys.}, 48(3):251--261, 1999.

\bibitem{MR}
P.~Mathonet and F.~Radoux.
\newblock Natural and projectively equivariant quantizations by means of
  {C}artan connections.
\newblock {\em Lett. Math. Phys.}, 72(3):183--196, 2005.

\bibitem{McDuff}
Dusa McDuff and Dietmar Salamon.
\newblock {\em Introduction to symplectic topology}.
\newblock Oxford Mathematical Monographs. The Clarendon Press Oxford University
  Press, New York, second edition, 1998.

\bibitem{OC}
V.~Ovsienko.
\newblock Vector fields in the presence of a contact structure.
\newblock {\em Enseign. Math.}, 52:215--229, 2006.

\bibitem{OvsBook0}
V.~Ovsienko and S.~Tabachnikov.
\newblock {\em Projective differential geometry old and new}, volume 165 of
  {\em Cambridge Tracts in Mathematics}.
\newblock Cambridge University Press, Cambridge, 2005.
\newblock From the Schwarzian derivative to the cohomology of diffeomorphism
  groups.

\bibitem{Wil}
E.~J. Wilczynski.
\newblock {\em Projective differential geometry of curves and ruled surfaces}.
\newblock Chelsea Publishing Co., New York, 1962.

\end{thebibliography}
\end{document}